\documentclass{amsart}

\usepackage{amssymb,color}
\usepackage{amsfonts}
\usepackage{amsmath}
\usepackage{euscript}
\usepackage{enumerate}
\usepackage{graphics}

\newtheorem*{Theorem1'}{Theorem 1'}

\theoremstyle{definition}

\theoremstyle{remark}

\numberwithin{equation}{section}

\newcommand \al{{\alpha}}

\begin{document}

\begin{center}
{\bf ROOTS MULTIPLICITY AND SQUARE-FREE FACTORIZATION OF POLYNOMIALS USING COMPANION MATRICES}
\end{center}

\medskip
\begin{center}{{\large {\sc N. H. Guersenzvaig}}\\
{\small A\!v. Corrientes 3985 6A, (1194) Buenos Aires, Argentina\\ email: nguersenz@fibertel.com.ar}}
\end{center}

\begin{center}
and
\end{center}

\begin{center}{{\large {\sc Fernando Szechtman}}\footnote{Corresponding author}\\
{\small Department of Mathematics and Statistics, University of Regina, Saskatchewan, Canada\\
email: fernando.szechtman@gmail.com}}
\end{center}

\vskip 0.3cm
\begin{center}{\large{ Abstract}}\end{center}

{\rm Given an arbitrary monic polynomial $f$ over a field $F$ of characteristic 0, we  use  companion matrices to construct a polynomial $M_f\in F[X]$ of minimum degree such that for each root $\alpha$ of $f$ in
the algebraic closure of~$F$, $M_f(\alpha)$ is equal to the multiplicity $m(\alpha)$ of $\alpha$ as a root of $f$.
As an application of $M_f$ we give a new method  to compute in $F[X]$ each component of the square-free factorization $f=P_1P_2^2\cdots P_m^m$, where $P_k$ is the product of all $X-\alpha$ with $m(\alpha)=k$, for $k=1, \dots , m=\max m(\alpha)$.




\section{Introduction}
\label{intr}


We fix throughout the paper a field $F$ of characteristic 0 and a monic polynomial $f\in F[X]$ of degree $n$. We say that $f$ is square-free if for any monic polynomial $g\in F[X]$: $g^2|f$ implies $g=1$, or, equivalently, if $\gcd(f, f')=1$. We assume henceforth that $n\geq 1$ and that $f$ has
 prime factorization
\begin{equation}
\label{canf}
f=f_{1}^{m_1}\cdots f_r^{m_r},
\end{equation}
so $f$ is square-free if and only if each $m_j=1$. Set $m=\max_{1\le j\le r} m_j$ and collect in $P_k$ all monic irreducible factors $f_j$ of $f$ having multiplicity $m_j=k$ in (\ref{canf}), that is,
$$P_k =\prod_{\substack{1\le j\le r\\ m_j=k}}f_j, \,\,k=1, \dots , m. $$
Each $P_k$ is square-free, yielding the so-called square-free factorization of $f$, namely
\begin{equation}
\label{sq}
f=P_{1}P^2_{2}\cdots P^m_{m}.
\end{equation}

Efficient algorithms to compute (\ref{sq}) were developed about forty years ago by Tobey \cite{7}, Horowitz \cite{4, 5}, Musser \cite{6} and Yun \cite{8}. A more recent algorithm can be found in \cite{1}. The square-free factorization  (\ref{sq}) is usually utilized as a first step towards the computation of the full factorization (\ref{canf}) as well as to
integrate rational functions, following a method due to Hermite \cite{3}, as explained in \cite{1}.

Let us review the Tobey-Horowitz algorithm, one of the simplest known methods to obtain (\ref{sq}).
Let $g=f_j$, where $1\leq j\leq r$, and set $k=m_j$ as well as
$h=f/g^k$. We then have $f=g^k h$ with $\gcd(g, h)=1$. Moreover,
$$
f'=kg^{k-1}g'h+g^k h'.
$$
As $g$ is relatively prime to $k g' h$, the multiplicity of $g$ in $\gcd(f, f')$ is $k-1$. Hence
\begin{equation*}
\label{eqq}
\gcd(f, f') =P_2P_3^2\cdots P_m^{m-1}.
\end{equation*}
Setting
$$
D_0=f,\; D_1=\gcd(D_0, D_0'),\; D_2=\gcd(D_1, D_1'),\dots
$$
a repeated application of the preceding discussion yields
$$D_k=\begin{cases}P_{k+1}P^2_{k+2}\cdots P^{m-k}_{m} &\text{$k=0, \dots , m-1$},\\
\,1  &\text{$k\geq m$,}
\end{cases}
$$
whence
\begin{equation}
\label{eqqq}
\frac{D_{k-1}}{D_{k}}=P_{k}P_{k+1}\cdots P_{m}, \, k=1, \dots , m-1.
\end{equation}
Therefore,  not only $m$ can be recognized as the smallest positive integer~$k$ such that $D_k=1$, but we also we have
$$
P_{k} = \frac{D_{k-1}}{D_{k}}:\frac{D_{k}}{D_{k+1}} =
\frac {D_{k-1}D_{k+1}}{D^2_k}, \,k=1, \dots , m.$$

In this paper we furnish an entirely new practical method to compute (\ref{sq}) by means of
Lagrange Interpolation Formula and companion matrices.

\smallskip
Let $K$ be an splitting field of $f$ over~$F$ and let $S(f)$ be the set of roots of $f$ in~$K$, whose
size will be denoted by~$s$.
For each $\al\in S(f)$ let $m(\al)$ be the multiplicity of~$\al$ as a root of $f$. Lagrange Interpolation Formula
ensures the existence and uniqueness of a polynomial
$M_f\in K[X]$ of degree less than $s$ satisfying
\begin{equation}
\label{eme1}
M_f(\alpha)=m(\alpha), \,\alpha\in S(f).
\end{equation}

The polynomial $f_0=f_1\cdots f_r$ is called the square-free part of $f$. Clearly,
$$
f_0=\underset{\alpha\in S(f)}\prod(X-\alpha)= P_1\cdots P_m.
$$
Applying case $k=1$ of (\ref{eqqq}) we obtain
$$
f_0= f/\gcd(f, f').
$$
Combining the two preceeding equations we obtain
 \begin{equation}
\label{ttr}
P_k = \prod_{\substack{\alpha\in S(f)\\ m(\alpha)=k}}(X-\alpha)= \gcd(M_f-k, \, f/\gcd(f, f')), \,\,k=1, \dots , m.
\end{equation}
The practical use of (\ref{ttr})  seems to be limited due to the mysterious nature of $M_f$. A closer look reveals that
\begin{equation}
\label{eme3}
M_f= \sum_{\alpha\in S(f)}\frac{m(\alpha)}{f'_0(\alpha)}\frac{f_0}{X-\alpha},
\end{equation}
which can be verified by evaluating both sides at each $\al\in S(f)$.
This seems to confirm the ineffectiveness of (\ref{ttr}), as (\ref{eme3}) requires the use of all roots of $f$ in~$K$.
However, applying the Galois group $G=\mathrm{Gal}(K/F)$ to (\ref{eme1}) we see that $M_f\in F[X]$. Indeed, let $\sigma\in G$.
Then $$M_f^\sigma(\alpha^\sigma)=M_f(\alpha)^\sigma=m(\alpha)^\sigma=m(\alpha)=m(\alpha^\sigma)=M_f(\alpha^\sigma),$$ so $M_f^\sigma=M_f$ by uniqueness, whence $M_f\in F[X]$, as claimed.
The fact that $M_f\in F[X]$ suggests that there should be a rational procedure to obtain $M_f$ from~$f$. This is exactly what we do in this paper, by means of companion matrices.

The outcome is an efficient procedure to compute all components $P_k$ of (\ref{sq}). Indeed, first compute $f_0=f/\gcd(f, f')$, then $M_f$,
as indicated in \S\ref{n2}, and then proceed sequentially to find all $P_1,P_2,\dots$ making use of (\ref{ttr}), stopping at $m$, namely
the smallest positive integer $k$ satisfying
$\deg(P_1)+ 2\deg(P_2)+\cdots  + k\deg(P_k)=n$.

\section{Using companion matrices to find $M_f$ from $f$}
\label{n2}

\vskip 0.1cm
 Given a monic polynomial $g=g_0+g_1X+\cdots+g_{s-1}X^{s-1}+X^{s}\in F[X]$ of positive degree $s$, its companion matrix $C_g\in M_s(F)$  is defined
$$
C_g=\left(%
\begin{array}{ccccc}
  0 & 0 & \cdots & 0 & -g_0 \\
  1 & 0 & \cdots & 0 & -g_1 \\
  0 & 1 & \cdots & 0 & -g_2 \\
  \vdots & \vdots & \cdots & \vdots & \vdots \\
  0 & 0 & \cdots & 1 & -g_{s-1} \\
\end{array}%
\right).
$$
We will use below the well-known fact that $g(C_g)=0$. Let $F_s[X]$ be the subspace of $F[X]$ with basis
$1, X,\dots, X^{s-1}$ and let $[R]\in F^s$ stand for
the coordinates of $R\in F_s[X]$ relative to this basis. We will also use the following formula from \cite{2}:
\begin{equation}
\label{form}
R(C_g)= ([R] \,\,C_g[R] \,\,C^2_g[R] \,\cdots \, C^{s-1}_g[R]).
\end{equation}
Since the first column of a product of matrices, say $AB$, is equal to the product of $A$ by the first column of $B$, it follows from (\ref{form})
that for any $P, Q, T\in F_s[X]$,
\begin{equation}
\label{x}
P(C_g)Q(C_g)= T(C_g) \iff P(C_g)[Q] = [T].
\end{equation}

\medskip
\noindent{\bf Theorem 2.1.} {\it Let $F$ be a field of characteristic 0 and let $f\in F[X]$ be monic of degree $n\geq 1$.
Set $f_0=f/\gcd(f, f')$, $s=\mathrm{deg}(f_0)$ and $P=f'/\gcd(f,f')$. Let $M_f$, $S(f)$ and $m(\alpha)$ be as in Section \ref{intr}.
Then there exist unique $g\in F_s[X]$ and $h\in F[X]$ such that $f'_0g + f_0h=1$. Moreover, for such $P$, $f_0$ and $g$,
$$[M_f]= P(C_{f_0})[g].$$}
\noindent{\sl Proof.} Definition of $f_0$ given in Section 1 guarantees $\gcd(f_0,f_0')=1$, which ensures the existence and uniqueness of $g$ and $h$.
Let $R$ be the remainder of dividing $M_f f'_0$ by $f_0$. By Lagrange Interpolation Formula, we have
$$
R =
 \sum_{\alpha\in S(f)}\frac{R(\alpha)}{f'_0(\alpha)}\frac{f_0}{X-\alpha}=
  f_0\sum_{\alpha\in S(f)}\frac{m(\alpha)}{X-\alpha}
  = f_0\frac{f'}{f}
  =P.
$$
Therefore,
$$M_{f}(C_{f_0})f'_0(C_{f_0}) = R(C_{f_0})=P(C_{f_0}).$$
But $f'_0(C_{f_0})g(C_{f_0})=I_s$, so
$$M_{f}(C_{f_0}) =P(C_{f_0})g(C_{f_0}),$$
which by (\ref{x}) is equivalent to what we want to prove. \qquad $\blacksquare$

\medskip
\noindent{\bf Example 2.2.} Let $f= X^4-4X + 3\in \mathbb Q[X]$. Then $$f'=4X^3-4\text{ and }\gcd(f,\, f')= X-1,$$
so
$$f_0=X^3+X^2+X-3, \,\,\,f'_0= 3X^2+2X +1 \text{ \,and }\,P=4X^2+4X+4.$$
We achieve $f'_0g + f_0h=1$ by taking
$$g =(1/72)X^2+(1/9)X +(1/24), \,\,\,h=(-1/24)X -(23/72).$$
Using
\begin{equation}
\label{x1}
C_{f_0}=\left(\begin{array}{ccc}
  0 & 0 & 3 \\
  1 & 0 & -1 \\
  0 & 1 & -1
\end{array}\right),
\end{equation}
we obtain
$$[M_f]=([P]\,\, C_{f_0}[P]\,\,C^2_{f_0}[P])[g]=
\left(\begin{array}{ccc}
  4 & 12 & 0 \\
  4 & 0 & 12 \\
  4 & 0 & 0
\end{array}
\right)\left(\begin{array}{c}
1/24 \\
1/9\\
1/72
\end{array}\right)= \left(\begin{array}{c}
 3/2 \\
1/3\\
1/6
\end{array}\right),$$
which means
\begin{equation}
\label{x2}
M_f= \frac{1}{6}X^2 + \frac 13X + \frac 32.
\end{equation}
Next we compute
$$P_{1}= \gcd(M_f-1, f_0)=  X^2+2X+3, \,\,\,P_{2}= \gcd(M_f-2, f_0) = X-1,$$
which yields the square-free factorization, $f=P_1P_2^2$.

\medskip

\noindent{\bf Remark 2.3.} It is possible to find the degrees of the $P_k$'s before actually computing these polynomials. This can be done as follows. Suppose $f_0$ has roots $\alpha_1, \dots , \alpha_s$ in $K$. Then $C_{f_0}$ is similar to the diagonal matrix \text{Diag}$(\alpha_1, \dots , \alpha_s)$, so $M_{f}(C_{f_0})$ is similar to \text{Diag}$(m(\alpha_1), \dots , m(\alpha_s))$ and has
characteristic polynomial
$$\text{Char}_{{M_f}(C_{f_0})}=\prod_{1\le i\le s} (X-m(\alpha_i))=\prod_{1\le k\le \max m(\alpha_i)}(X-k)^{\deg(P_{k})}.$$
For instance, in Example 2.2 we can use (\ref{x1}) and (\ref{x2}) to obtain
$$
M_f(C_{f_0})=\left(\begin{array}{ccc}
3/2 & 1/2 & 1/2 \\
1/3 & 4/3 & 1/3 \\
1/6 & 1/6 & 7/6
\end{array}
\right).
$$
Therefore
$$\text{Char}_{M_f(C_{f_0})}= X^3-4X^2+5X-2 = (X-1)^2(X-2),$$
which indicates the square-free factorization $f=P_1P_2^2$, with $\mathrm{deg}(P_1)=2$ and $\mathrm{deg}(P_2)=1$,
in complete agreement with what we found above.

\medskip
\noindent{\bf Remark 2.4.} Let $R$ be an integrally closed domain of characteristic 0 with field of fractions $F$, and let $f\in R[X]$ be monic of degree $n\geq 1$. Then the square-free components of $f$, as found above, will have coefficients in $R$. Indeed, (\ref{sq})
shows that all $\alpha_i$, and hence the coefficients of all $P_k$, are integral over $R$. But these coefficients are in $F$, and hence lie in $R$.



\medskip

\noindent{\bf Acknowledgements.}  We thank V. Sergeichuk for valuable comments.


\end{document}